\documentclass[english]{article}

\usepackage{babel}
\usepackage{amstext}
\usepackage{amsmath}
\usepackage{amsfonts}
\usepackage{latexsym}
\usepackage{ifthen}
\usepackage{xypic}
\xyoption{all}
\newcommand{\id}{{\rm id}}
\newcommand{\ball}{{\mathbb B}}
\renewcommand{\O}{{\cal O}}
\newcommand{\PB}{{\mathbb P}}

\newcommand{\Pic}{{\rm Pic}}
\newcommand{\GHom}{{\cal H}om}
\newcommand{\Hom}{{\rm Hom}}
\newcommand{\Ext}{{\rm Ext}}
\newcommand{\End}{{\cal E}nd}
\newcommand{\PN}[1]{{\mathbb P}_{#1}}
\newcommand{\lra}{\longrightarrow}
\newcommand{\KC}{{\mathbb C}}

\newcommand{\Formel}[1]{(\ref{#1})}
\newtheorem{lemma1}[equation]{}

\newenvironment{example}{\begin{lemma1}{\bf Example.}\rm}{\end{lemma1}}
\newenvironment{abs}{\begin{lemma1}\rm}{\end{lemma1}}
\newenvironment{theorem}{\begin{lemma1}{\bf Theorem.}}{\end{lemma1}}
\newenvironment{theorem2}[1]{\begin{lemma1}{\bf Theorem [#1].}}{\end{lemma1}}
\newenvironment{proposition}{\begin{lemma1}{\bf Proposition.}}{\end{lemma1}}
\newenvironment{corollary}{\begin{lemma1}{\bf Corollary.}}{\end{lemma1}}
\newenvironment{remark}{\begin{lemma1}{\bf Remark.}\rm}{\end{lemma1}}

\newenvironment{proof}{\noindent {\em Proof}.}{\hfill $\Box$}
\newenvironment{proof2}[1]{\noindent {\em Proof of #1}.}{\hfill $\Box$}

\newcommand{\Abs}[1]{\ref{#1}}
\newcommand{\Theo}[1]{Theorem~\ref{#1}}
\newcommand{\Prop}[1]{Proposition~\ref{#1}}
\newcommand{\Cor}[1]{Corollary~\ref{#1}}

\begin{document}

\title{Splitting jet sequences}
\author{Priska Jahnke and Ivo Radloff\thanks{The authors were supported by a Forschungsstipendium of the Deutsche Forschungsgemeinschaft.}}
\date{October 2002}
\maketitle

\section*{Introduction}

It is well known, that the splitting of the first jet sequence of a vector bundle $E$ on a complex manifold $X$ implies that $E$ admits a holomorphic affine connection. In this note we show that the splitting of a higher jet sequence 
  \[0 \lra S^k\Omega_X^1\otimes E \lra J_k(E) \lra J_{k-1}(E) \lra 0,\]
where $k > 1$, does not only have consequences for $E$, but also has strong consequences for the base manifold $X$. The following result is in a sense an infinitesimal converse to a result of Biswas, \cite{Bis99}, saying that if $X$ is a complex manifold admitting a holomorphic projective structure, then the higher jet sequences of certain tensor powers of a theta characteristic on $X$ split holomorphically in a charateristic manner (see \Abs{psjets}):

\vspace{0.2cm}

\noindent {\bf Theorem.} {\em Let $X$ be a compact K\"ahler manifold of dimension $n$ and $E$ a vector bundle of rank $r$ on $X$. If the $k$--th jet sequence of $E$ splits for some $k > 1$, then $X$ admits a holomorphic normal projective connection and $J_{k-1}(E)$ admits a holomorphic affine connection. Moreover
 \begin{enumerate}
  \item if $\det E^*$ is not nef and $X$ is projective or if $X$ contains a rational curve, then $X \simeq \PN{n}$ and $E \simeq \oplus \O_{\PN{n}}(k-1)$,
  \item if $\det E^* \equiv 0$, then $X$ is covered by a torus and $E$ admits a holomorphic affine connection,
  \item if $\det E^*$ is ample, then $X$ is covered by the unit ball.
 \end{enumerate} 
Here $E^* = \GHom(E, \O_X)$ denotes the dual of $E$.}

\vspace{0.2cm}

\noindent The Theorem implies for example that the Chern classes of a compact K\"ahler manifold $X$ carrying a vector bundle $E$ with splitting $k$--th jet sequence for some $k >1$ are related as follows:
 \[(n+1)^{\mu} c_{\mu}(E) = \mbox{$\left({r \atop {\mu}}\right)$} (k-1)^{\mu} c_1(T_X)^{\mu}, \quad 
    r^{\mu}(k-1)^{\mu} c_{\mu}(T_X) = \mbox{$\left({n+1 \atop \mu}\right)$} c_1(E)^{\mu}.\]

Motivation to this note was the following problem. Projective space $\PN{n}$, \'etale quotients of complex tori and manifolds, whose universal cover is the unit ball in $\KC^n$ are basic examples of complex compact K\"ahler manifolds carrying a flat holomorphic normal projective connection. By \cite{PC} this is not a complete list. Whether any compact complex K\"ahler manifold with a holomorphic normal projective connection carries a holomorphic projective structure, as in the affine case, seems still to be an open problem. While the uniformization theorem says that the answer to the problem is ``yes'' in the case of curves, Kobayashi and Ochiai showed that it is as well ``yes'' in the case of surfaces (\cite{KoOc}). It is as well ``yes'' for projective threefolds (\cite{PC}). 

Splitting jet sequences may be a key to the understanding of this problem. In the affine case, we find for example:

\vspace{0.2cm}

\noindent {\bf Corollary.} {\em Let $X$ be a compact K\"ahler manifold. The first jet sequence of $\O_X$ is always split exact. One of the higher jet sequences
 \begin{equation} \label{jetOX}
   0 \lra S^k\Omega^1_X \lra J_k(\O_X) \lra J_{k-1}(\O_X) \lra 0,
 \end{equation}
$k>1$, splits if and only if $X$ is covered by a torus, i.e., admits a flat holomorphic affine connection. The splitting of one higher jet sequence of $\O_X$ implies the splitting of all jet sequences of $\O_X$.}
\vspace{0.2cm}

\noindent The question is whether a similar result holds if one replaces $\O_X$ in \Formel{jetOX} by certain powers of a (numerical) theta characteristic on $X$ as suggested by Biswas's result.


\section{Projective structures and connections}
\setcounter{equation}{0}

In this section we briefly recall some basic definitions and results. Throughout this section, we assume that $X$ is a compact K\"ahler manifold of dimension $n$, even though the K\"ahler condition is not always necessary. Concerning our notations see \Abs{notation}.

\begin{abs} \label{pstructure}
{\bf Holomorphic affine and projective structures.} (\cite{KoOc} or \cite{KoWu}). $X$ is said to admit a holomorphic {\em affine} (resp. {\em projective}) structure, if it can be covered by coordinate charts, such that the coordinate change is given by restrictions of holomorphic affine (resp. projective) transformations of $\KC^n$ (resp. $\PN{n}$). A manifold with a holomorphic affine structure admits a holomorphic projective structure. The ``list of standard examples'' of compact complex K\"ahler (--Einstein) manifolds admitting a holomorphic projective structure, already mentioned above, consists of
\begin{enumerate}
 \item the projective space $\PN{n}$,
 \item \'etale quotients of complex tori,
 \item manifolds, whose universal cover is the unit ball ${\ball}^n \subset \KC^n$.
\end{enumerate}
In 3.) note that the group of automorphisms of ${\ball}^n$ is ${\rm SU}(1,n)$, a subgroup of ${\rm PGl}(n+1)$. By the uniformization theorem, every compact complex curve is in the list of standard examples, i.e., admits a holomorphic projective structure. The manifolds in 2.) admit a holomorphic affine structure.

There is one more (projective) example in dimension $3$, namely {\'e}tale quotients of smooth modular families of false elliptic curves parametrized by a Shimura curve (\cite{PC}).
\end{abs}

The corresponding infinitesimal notions to holomorphic affine and projective structures are holomorphic affine and projective connections. Before we come to this we recall the definition of the Atiyah class of a holomorphic vector bundle:

\begin{abs} \label{Atiyahclass}
{\bf The Atiyah class and $a(E)$.} (\cite{Ati}). Let $E$ be a holomorphic vector bundle of rank $r$ on $X$ and $\{U_{\alpha}; z_{\alpha}^1, \dots, z_{\alpha}^n\}_{\alpha \in I}$ coordinate charts with coordinates $z_{\alpha}^i$ where $E$ is trivial. Let $\{U_{\alpha}; e_{\alpha}^1, \dots, e_{\alpha}^r\}_{\alpha \in I}$ be a local frame for $E$. Let $g_{\alpha\beta} \in H^0(U_{\alpha} \cap U_{\beta}, {\rm Gl}(r, \O_X))$ be transition functions of $E$ such that $e_{\beta}^k = \sum_l g_{\alpha\beta}^{lk}e_{\alpha}^l$.

The {\em Atiyah class of $E$} is the splitting obstruction of the first jet sequence, i.e., it is the image of $\id_E$ under the first connecting morphism 
 \[H^0(X, \End(E)) \lra \Ext^1(E, \Omega_X^1 \otimes E) \simeq H^1(X, \End(E) \otimes \Omega^1_X).\] 
The Dolbeault isomorphism $H^1(X, \End(E) \otimes \Omega_X^1) \simeq H^{1,1}(X, \End(E))$ maps the Atiyah class to $[-\Theta_h]$, where $\Theta_h$ denotes the canonical curvature of $E$ with respect to a hermitian metric $h$ on $E$. In particular, the trace of the Atiyah class is $-2\pi i c_1(E)$ in $H^1(X, \Omega^1_X)$.

If we define $a(E)$ as $-\frac{1}{2\pi i}$ times the Atiyah class of $E$, then the trace of $a(E)$ is $c_1(E)$, which makes this definition convenient for our purposes. In local coordinates, $a(E)$ is the class of the Chech cocycle $a(E)_{\alpha\beta} \in Z^1({\cal U}, \End(E)\otimes \Omega_X^1)$, where
 \[a(E)_{\alpha\beta} = \mbox{$\frac{1}{2\pi i} \sum_{i,j,l} \frac{\partial g_{\alpha\beta}^{jl}}{\partial z_{\alpha}^i}dz_{\alpha}^i \otimes e_{\alpha}^j \otimes e_{\beta}^{l*} = \frac{1}{2\pi i} \sum_{1 \le j,l \le r} dg_{\alpha\beta}^{jl} e_{\alpha}^j \otimes e_{\beta}^{l*}$}.\]
See \cite{Ati} for the functorial behavior of $a(E)$ under pull--back, tensor products and direct sums.
\end{abs}

\begin{abs} \label{pc}
{\bf Holomorphic affine and projective connections.} (\cite{KoOc} or \cite{KoWu}) $X$ is said to have a {\em holomorphic affine connection} if
 \[a(\Omega^1_X) = 0 \;\mbox{ in } \; H^1(X, \Omega_X^1 \otimes T_X \otimes \Omega_X^1),\]
where $T_X$ and $\Omega_X^1$ denote the holomorphic tangent sheaf and the sheaf of holomorphic $1$--forms, respectively. If $X$ has a holomorphic affine structure, then $X$ has a (flat) holomorphic affine connection. Since we assume $X$ K\"ahler, the existence of a holomorphic affine connection implies that $X$ is covered by a torus (\cite{KoWu}, 2.4.1.~Theorem). In other words: a compact K\"ahler manifold admits a holomorphic affine connection if and only if it is covered by a torus.

$X$ is said to have a {\em holomorphic (normal) projective connection} if
 \[a(\Omega_X^1) =  \mbox{$\id_{\Omega_X^1} \otimes \frac{c_1(K_X)}{n+1} + \frac{c_1(K_X)}{n+1} \otimes \id_{\Omega_X^1}$} \; \mbox{ in } \; H^1(X, \Omega_X^1 \otimes T_X \otimes \Omega_X^1),\]
where we use the identifications
 \[\End(\Omega^1_X) \otimes \Omega_X^1 \simeq \Omega_X^1 \otimes T_X \otimes \Omega_X^1 \simeq \Omega_X^1 \otimes \End(\Omega^1_X)\]
and consider $c_1(K_X)$ as an element in $H^1(X, \Omega_X^1)$. (See \cite{MoMo} for a more differential geometric description of projective connections). From now on we will drop the appellation ``normal''. If $X$ has a holomorphic projective connection and $c_1(T_X) = 0$ in $H^1(X, \Omega_X^1)$, then $X$ has a holomorphic affine connection. Hence: a compact K\"ahler manifold with a holomorphic projective connection and $c_1(T_X) = 0$ is covered by a torus.

If $X$ has a holomorphic affine (projective) structure, then $X$ admits a holomorphic affine (projective) connection. A holomorphic affine (projective) connection is said to be {\em flat} or {\em integrable} if it corresponds to a holomorphic affine (projective) structure. The examples in \Abs{pstructure} all carry a flat holomorphic projective connection.

Gunning's formula on the Chern classes of a K\"ahler manifold with a holomorphic projective connection says (\cite{Gun}, p.94)
  \begin{equation} \label{ChernClss}
    (n+1)^{\mu}c_{\mu}(T_X) = \mbox{$\left({n+1 \atop \mu}\right)$} c_1^{\mu}(T_X) \; \mbox{ in } \; H^{\mu}(X, \Omega_X^{\mu}).
  \end{equation}
Important for our purposes will be moreover the following result of Kobayashi and Ochiai:

\begin{theorem2}{Kobayashi, Ochiai} \label{KE}
 The list of K\"ahler--Einstein manifolds admitting a holomorphic projective connection is the list of standard examples.
\end{theorem2}

In {\em Holomorphic projective structures and compact complex surfaces I, II}, \cite{KoOc}, Kobayashi and Ochiai show that the list of compact K\"ahler surfaces admitting a holomorphic projective connection is precisely the list of standard examples. It was proved in {\em On manifolds with holomorphic normal projective connections}, \cite{PC}, that even in the case of projective threefolds the notion of holomorphic projective structures and connections coincide and that there is besides the standard examples only one more type that completes the list: {\'e}tale quotients of smooth modular families of false elliptic curves.
\end{abs}

\begin{abs} \label{psjets}
{\bf Projective structures and splitting jet sequences.} In {\em Differential operators on complex manifolds with a flat projective structure} (\cite{Bis99}) I.~Biswas shows, that if $X$ is a complex manifold equipped with a projective structure and $\theta$ is a theta characteristic on $X$, i.e. a  holomorphic line bundle such that $\theta^{n+1} \simeq K_X$, then, for any $k \ge 0$, the jet bundle $J_k(\theta^{-k})$ has a natural flat connection and the restriction homomorphism,
 \[J_l(\theta^{-k}) \lra J_k(\theta^{-k})\]
admits a canonical splitting for $l \ge k$ (\cite{Bis99}, Theorem~3.7.). In particular, the $(k+1)$--th jet sequence of $\theta^{-k}$ splits in the situation of the theorem.
\end{abs}

\begin{abs} \label{notation}
\noindent {\bf Notation and conventions.} For a compact complex manifold $X$, the canonical divisor is denoted by $K_X$. We will identify line bundles and divisors and write $K_X$ instead of $\O_X(K_X)$. The tensor product of line bundles will be denoted by $+$ or $\otimes$. A line bundle is called nef, if the intersection number with every irreducible curve is non--negative. It is called big, if the top self intersection class is positive. A vector bundle is called nef, if $\O_{\PB(E)}(1)$ is nef on $\PB(E)$. The symbol $\equiv$ denotes numerical equivalence. A theta characteristic is a line bundle $\theta \in \Pic(X)$, such that $\theta^{n+1} = K_X$. 
\end{abs}


\section{Splitting jet sequences} \label{sec jets}
\setcounter{equation}{0}

The main result of this section is the following first order consequence of the splitting of jet sequences:

\begin{theorem} \label{AtClass}
   Let $E$ be a vector bundle on the compact complex manifold $X$. If the $k$--th jet sequence of $E$ splits for some $k \ge 1$, then 
 \begin{enumerate}
  \item the $(k-1)$--th jet bundle of $E$ admits a holomorphic affine connection, i.e., $a(J_{k-1}(E)) = 0$ in $H^1(X, \End(J_{k-1}(E)) \otimes \Omega_X^1)$, 
 \item in $H^1(X, \End(E) \otimes \Omega_X^1 \otimes T_X \otimes \Omega_X^1)$ holds
  \[(k-1) \cdot \id_E \otimes a(T_X)  = a(E) \otimes \id_{T_X} + \id_{T_X} \otimes a(E).\]
 \end{enumerate}
\end{theorem}

\begin{proof}
Let $a_k(E)$ be $-\frac{1}{2i\pi}$ times the splitting obstruction of the $k$--th jet sequence of $E$. Then $a_k(E)$ is an element in $\Ext^1(J_{k-1}(E), S^k\Omega^1_X \otimes E)$ and $a_1(E) = a(E)$ from \Abs{Atiyahclass}.

\vspace{0.2cm}

\noindent 1.) The first jet sequence of $J_{k-1}(E)$ and the $k$--th of $E$ give the diagram:
  \[\xymatrix{0 \ar[r] & S^k\Omega^1_X \otimes E   \ar@{^{(}->}[d] \ar[r] & J_k(E) \ar[r] \ar@{^{(}->}[d] & J_{k-1}(E) \ar@{=}[d]\ar[r] & 0 \\
 0 \ar[r] & \Omega^1_X \otimes J_{k-1}(E) \ar[r] & J_1(J_{k-1}(E)) \ar[r] & J_{k-1}(E) \ar[r] & 0 \\}\]
Applying $\Hom(J_{k-1}(E), -)$ to the diagram gives a map
  \[\Ext^1(J_{k-1}(E), S^k\Omega^1_X \otimes E) \lra \Ext^1(J_{k-1}(E), \Omega^1_X \otimes J_{k-1}(E)),\]
mapping $a_k(E)$ to $a_1(J_{k-1}(E)) = a(J_{k-1}(E))$. This shows the first part of the theorem.

\vspace{0.2cm}

\noindent 2.) If $k = 1$, i.e., if the first jet sequence of $E$ splits, then $a(E) = a(J_0(E)) = 0$ by 1.), and the claim is trivial. So assume $k > 1$. Applying $\Hom(-, S^k\Omega_X^1 \otimes E)$ to the $(k-1)$--th jet sequence of $E$ gives a map
  \begin{equation} \label{abar}
   \Ext^1(J_{k-1}(E), S^k\Omega_X^1 \otimes E) \lra \Ext^1(S^{k-1}\Omega_X^1 \otimes E, S^k\Omega_X^1 \otimes E).
  \end{equation}
The image of $a_k(E)$ under this map will be called $\bar{a}_k(E)$. Since the $k$--th jet sequence of $E$ splits by assumption, $\bar{a}_k(E) = 0$. A contraction map
 \[\Ext^1(S^{k-1}\Omega_X^1 \otimes E, S^k\Omega_X^1 \otimes E) \lra H^1(X, \End(E) \otimes \Omega_X^1 \otimes T_X \otimes \Omega_X^1)\]
will show the claim.

Compute now $\bar{a}_k(E)$ locally: choose coordinate charts $\{U_{\alpha}; z_{\alpha}^1, \dots, z_{\alpha}^n\}_{\alpha \in I}$ in $X$, such that $E$ is trivial on each $U_{\alpha}$ with local frame $\{U_{\alpha}; e_{\alpha}^1, \dots, e_{\alpha}^r\}_{\alpha \in I}$ and transition functions $g_{\alpha\beta}$. Let $h_{\alpha\beta}^i(z_{\alpha}^1, \dots, z_{\alpha}^n) = z_{\beta}^i$ on $U_{\alpha\beta}$ be a holomorphic coordinate change. Then 
 \[\left\{U_{\alpha}; \{dz_{\alpha}^{i_1} \cdots dz_{\alpha}^{i_k} \otimes e_{\alpha}^l\}_{{i_1 \le \dots \le i_k \atop 1 \le l \le r}}, \dots, \{dz_{\alpha}^i \otimes e_{\alpha}^l\}_{{1 \le i \le n \atop 1 \le l \le r}}, \{e_{\alpha}^l\}_l\right\}_{\alpha \in I}\]
is a local frame for $J_k(E)$. Denote by $j_kg_{\alpha\beta}$ the transition functions of $J_k(E)$. Then
 \[j_kg_{\alpha\beta} = \left(\begin{array}{cc}
     S^k \frac{\partial z_{\beta}}{\partial z_{\alpha}} \otimes g_{\alpha\beta} & v_{\alpha\beta}\\
      0 & j_{k-1}g_{\alpha\beta}
  \end{array}\right),\]
where $v_{\alpha\beta}: J_{k-1}|_{U_{\beta}} \to S^k\Omega_X^1 \otimes E|_{U_{\alpha}}$ defines the extension class $a_k(E)$. We are only interested in the image $\bar{a}_k$, defined in \Formel{abar}. This means that we will apply $v_{\alpha\beta}$ not to sections of $J_{k-1}(E)$ but only to sections of the subbundle $S^{k-1}\Omega_X^1 \otimes E$ on $U_{\beta}$. Denote the corresponding columns of $v_{\alpha\beta}$ by $\bar{v}_{\alpha\beta}$. Then 
 \[\bar{a}_k = [\bar{v}_{\alpha\beta}] \in \Ext^1(S^{k-1}\Omega_X^1 \otimes E, S^k\Omega_X^1 \otimes E).\]

The map $\bar{v}_{\alpha\beta}$ can be computed directly in terms of the local frames of $\Omega_X^1$ and $E$ on $U_{\alpha}$ and $U_{\beta}$. For example in the case $k = 2$:
 \[\bar{v}_{\alpha\beta} = \id_E \otimes \sum_{{i \le j \atop u}} \frac{\partial^2 z_{\beta}^u}{\partial z_{\alpha}^i \partial z_{\alpha}^j} dz_{\alpha}^i dz_{\alpha}^j \otimes dz_{\beta}^{u*} + \!\!\!\!\!\sum_{{{i \le j \atop \{r,s\} = \{i,j\}} \atop l,m}} \!\!\!\!\! \frac{\partial g_{\alpha\beta}^{lm}}{\partial z_{\alpha}^r} e_{\alpha}^l \otimes e_{\beta}^{m*} \otimes dz_{\alpha}^i dz_{\alpha}^j \otimes dz_{\alpha}^{s*}.\]
In the general case the first sum runs over all $i_1 \le \dots \le i_k$ and all pairs $(i_r, i_s)$, $r\not= s$ and all $u$, and each $dz_{\alpha}^i dz_{\alpha}^j \otimes dz_{\beta}^{u*}$ has to be replaced by 
 \[dz_{\alpha}^{i_1} \!\!\cdots dz_{\alpha}^{i_k} \otimes dz_{\beta}^{u*} dz_{\alpha}^{i_1*} \!\!\!\cdots \widehat{dz_{\alpha}^{i_r*}} \!\!\!\cdots \widehat{dz_{\alpha}^{i_s*}} \!\!\!\cdots dz_{\alpha}^{i_k*}.\] 
The second sum has to be modified analogously.

Let $V$ be a vector space of dimension $r$ with basis $b_1, \dots, b_r$. The natural maps
  \[A_k: S^{k-1}V^* \otimes S^{k-2}V \lra S^kV^* \otimes S^{k-1}V\] 
defined by $b_{i_1}^* \cdots b_{i_{k-1}}^* \otimes b_{u_1} \cdots b_{u_{k-2}} \mapsto \sum_{s=1}^r b_{i_1}^* \cdots b_{i_{k-1}}^*b_s^* \otimes b_{u_1} \cdots b_{u_{k-2}}b_s$ split; the splitting map $C_k: S^kV^* \otimes S^{k-1}V \to S^{k-1}V^* \otimes S^{k-2}V$ can be expressed as a linear combination of contraction maps. Therefore all composition maps
 \[F_{k,i} = A_k \circ \cdots \circ A_i : S^iV^* \otimes S^{i-1}V \lra S^kV^* \otimes S^{k-1}V, \quad i < k,\] 
admit a natural splitting, namely $S_{k,i} = C_i \circ \cdots \circ C_k$. 

Applied to the above situation, where $V$ is the vector space generated by $dz_{\alpha}^{i*}$, $i = 1, \dots, n$, we get
 \[\bar{v}_{\alpha\beta} = (k-2)!^{-1}F_{k,2}(\varphi_2) + (k-1)!^{-1}F_{k,1}(\varphi_1),\]
where $\varphi_2 \in S^2\Omega_X^1|_{U_{\alpha}}\otimes T_X|_{U_{\beta}}$ and $\varphi_1 \in \Omega_X^1|_{U_{\alpha}}$ are given by
 \[\varphi_2 = \id_E \otimes \sum_{u, i, j} \frac{\partial^2 z_{\beta}^u}{\partial z_{\alpha}^{i} \partial z_{\alpha}^{j}} dz_{\alpha}^i dz_{\alpha}^j \otimes dz_{\beta}^{u*},  \quad
   \varphi_1 = \sum_{i,j,l} \frac{\partial g_{\alpha\beta}^{jl}}{\partial z_{\alpha}^i}dz_{\alpha}^i \otimes e_{\alpha}^j \otimes e_{\beta}^{l*}.\]
The map $T_2: S^2V^*\otimes V \to V^* \otimes V^* \otimes V$, defined by mapping $b_ib_j \otimes b_l^*$ to $\frac{1}{2}(b_i \otimes b_j + b_j \otimes b_i) \otimes b_l^*$, gives
 \[2(k-1)!T_2S_{k,2}\bar{v}_{\alpha\beta} = (k-1) \id_E \otimes a(\Omega_X^1)_{\alpha\beta} + \id_{T_X} \otimes a(E)_{\alpha\beta} + a(E)_{\alpha\beta} \otimes \id_{T_X}.\]
If the $k$--th jet sequence of of $E$ splits, then $[\bar{v}_{\alpha\beta}] = 0$, i.e., the class of the left hand side is zero. Since $a(\Omega_X^1) = -a(T_X)$, the claim follows.
\end{proof}

\begin{corollary} \label{ProjConn}
  Let $E$ be a vector bundle on the compact K\"ahler manifold $X$. If the $k$--th jet sequence of $E$ splits for some $k > 1$, then $X$ admits a holomorphic projective connection and
  \begin{equation} \label{AtForm}
   (k-1) \cdot \id_E \otimes c_1(T_X)  =  (n+1) \cdot a(E) \quad \mbox{ in } H^1(X, \End(E) \otimes \Omega_X^1).
  \end{equation}
 In particular, if $X$ is a surface, a projective threefold or if $X$ admits a K\"ahler--Einstein metric, then $X$ admits a projective structure and is one of the manifolds from section~1.1.
\end{corollary}
The identity \Formel{AtForm} is obviously true even in the case $k = 1$, since the splitting of the first jet sequence means nothing but that $a(E) = 0$.

\begin{remark}
The corollary should be understood as follows. Assume that $X$ carries a theta characteristic $\theta \in \Pic(X)$. If $E$ is a vector bundle on $X$, such that the $k$--th jet sequence of $E$ splits for some $k > 1$, then the $i$--th jet sequence of $E \otimes \theta^{k-i}$ splits for $i = 1, \dots, k$. We will not prove this here; it will not be used in the sequel. \Cor{ProjConn} proves the case $i = 1$: the splitting of the first jet sequence of $E \otimes \theta^{k-1}$ means the vanishing of $a(E \otimes \theta^{k-1})$ or, equivalently, $a(E) = \id_E \otimes \frac{k-1}{n+1} \cdot c_1(T_X)$.
\end{remark}

\begin{proof2}{\Cor{ProjConn}}
 From the proof of the Theorem we know that the class of
 \[(k-1) \cdot \id_E \otimes a(\Omega_X^1)_{\alpha\beta} + \id_{T_X} \otimes a(E)_{\alpha\beta} + a(E)_{\alpha\beta} \otimes \id_{T_X}\]
in $H^1(X, \End(E) \otimes \Omega_X^1 \otimes T_X \otimes \Omega_X^1)$ is zero. For the identity \Formel{AtForm} consider the contraction of the last two factors in $\End(E) \otimes \Omega_X^1 \otimes T_X \otimes \Omega^1_X$. Then 
 \[(k-1) \cdot \id_E \otimes c_1(T_X) = n \cdot a(E) + a(E)\]
since $\id_{T_X} \otimes a(E)$ maps to $n \cdot a(E)$ and $a(E) \otimes \id_{T_X}$ to $a(E)$. This shows \Formel{AtForm}. Consider on the other hand the contraction of the first two factors $\End(E) \simeq E \otimes E^*$:
 \[r(k-1) \cdot a(T_X) =  \id_{T_X} \otimes c_1(E) + c_1(E) \otimes \id_{T_X}.\]
The claim follows now from the first set of equalities in the following \Cor{Chern}, which are a direct consequence of the identity \Formel{AtForm}. 
\end{proof2}

\

Using Gunning's formula \Formel{ChernClss} on the Chern classes of a compact K\"ahler manifold with projective connection and formula \Formel{AtForm} in \Cor{ProjConn} we obtain:

\begin{corollary} \label{Chern}
  In the situation of \Cor{ProjConn}, the Chern classes of $E$ and $X$ are related as follows:
 \[(n+1)^{\mu} c_{\mu}(E) = \mbox{$\left({r \atop \mu}\right)$} (k-1)^{\mu} c_1(T_X)^{\mu}, \quad 
    r^{\mu}(k-1)^{\mu} c_{\mu}(T_X) = \mbox{$\left({n+1 \atop \mu}\right)$} c_1(E)^{\mu}.\]
 For $k = 1$ this just means that the Chern classes of a vector bundle with trivial Atiyah class vanish.
\end{corollary}

\Cor{Chern} proves the first part of the theorem in the introduction. Before we come to the proof of the second part, we give the following

\begin{example}
 {\em If $E$ is an indecomposable vector bundle of rank $r$ on the smooth compact curve $C$, then the $k$--th jet sequence of $E$ splits if and only if
 \begin{equation} \label{curves}
   2 \cdot c_1(E) = r(k-1) \cdot c_1(T_C) = 2r(k-1)(1-g(C)).
 \end{equation}}
Note that if $E$ is decomposable, i.e., $E = \oplus E_i$ for some indecomposable vector bundles $E_i$, then the $k$--th jet sequence of $E$ splits if and only if the $k$--th jet sequence of $E_i$ splits for all $i$.

Indeed, if the $k$--th jet sequence of $E$ splits, then formula \Formel{curves} holds by \Cor{Chern}. Conversely, assume that $E$ is a vector bundle on $C$ satisfying $2 c_1(E) = r(k-1) \cdot c_1(T_C)$. Then $c_1(E \otimes \theta^{k-1}) = 0$, where $\theta$ denotes a theta characterstic on $C$. By \cite{Ati}, {\em Proposition}~14 and 19 the bundle $N = E \otimes \theta^{k-1}$ can then be written with constant transition functions. Therefore $J_k(E) = J_k(\theta^{1-k}) \otimes N$ for all $k$. The claim follows now from Biswas' result. 
\end{example}


\section{Classification}
\setcounter{equation}{0}

We can now either consider special manifolds $X$ carrying a holomorphic vector bundle $E$ where one of the higher jet sequences split and study $E$, or, conversely, we can consider special bundles and ask what this splitting means for $X$. We start with the first case. In \cite{PC}, Ye's result whereafter projective space is the only smooth Fano variety $X$ admitting a holomorphic projective connection (\cite{Ye}) was generalized to the case of a compact K\"ahler manifold $X$ that contains a rational curve. Using only direct consequences of the splitting of jets we prove here:

\begin{proposition} \label{RatCurve}
 Let $X$ be a compact K\"ahler manifold and $E$ a vector bundle of rank $r$ on $X$. If $X$ contains a rational curve and if the $k$--th jet sequence of $E$ splits for some $k > 1$, then $X \simeq \PN{n}$ and $E \simeq \oplus \O_{\PN{n}}(k-1)$.
\end{proposition}

\begin{remark}
 Every vector bundle on $\PN{1}$ splits into a sum of line bundles.
\end{remark}

\begin{proof2}{\Prop{RatCurve}}
  Let $f: \PN{1} \to X$ be any rational curve in $X$.  We want to show $X = \PN{n}$: by \Theo{AtClass}, $a(J_{k-1}(E)) = 0$, hence $a(f^*J_{k-1}(E)) = 0$ on $\PN{1}$. Then the vector bundle $f^*J_{k-1}(E)$ arises from a representation of the fundamental group of $\PN{1}$ by \cite{Ati} Theorem~8, which means $f^*J_{k-1}(E)$ is trivial. Consider the jet bundle diagrams from \cite{Gro}:
\[\xymatrix{0 \ar[r] & S^{k-1}\Omega_{\PN{1}}^1 \otimes f^*E \ar[r] & J_{k-1}(f^*E) \ar[r] & J_{k-2}(f^*E) \ar[r] & 0  \\
 0 \ar[r] & f^*S^{k-1}\Omega_X^1 \otimes f^*E \ar[r]\ar[u] & f^*J_{k-1}(E) \ar[r]\ar[u] & f^*J_{k-2}(E) \ar[r]\ar[u] & 0.}\]
Since $f^*S^{\mu}\Omega^1_X \otimes f^*E \to S^{\mu} \Omega^1_{\PN{1}} \otimes f^*E$ is generically surjective for all $\mu$, it follows by induction that all the vertical maps are generically surjective. Hence  
\[0 \lra J_{k-1}(f^*E)^* \lra f^*J_{k-1}(E)^* = \O_{\PN{1}}^{\oplus m}\]
is injective. Therefore $J_{k-1}(f^*E)$ is nef, meaning that if we write $J_{k-1}(f^*E)$ as a sum of line bundles, non of them will have negative degree. Write $f^*E$ as a sum of line bundles $\O_{\PN{1}}(a_i)$. Since $J_{k-1}(f^*E)$ is nef, $J_{k-1}(\O_{\PN{1}}(a_i))$ is nef for all $i$. A direct computation shows: if $J_{k-1}(\O_{\PN{1}}(a))$ is nef on $\PN{1}$ for some $k > 1$, then $a > 0$. Hence $f^*E$ is ample. 
Dualizing the pull back of the $(k-1)$--th jet sequence of $E$ yields a surjective map
 \[f^*J_{k-1}(E)^* \otimes f^*E = f^*E^{\oplus m} \lra f^*(S^{k-1}T_X \otimes \End(E)).\]
Since $f^*E$ is ample and $\O_{\PN{1}}$ is a direct summand of $f^*\End(E)$, $f^*S^{k-1}T_X$ is ample. Since $k > 1$, $f^*T_X$ is ample. It follows now that $X$ is projective and by Mori's Theorem (\cite{Mori} and \cite{MiPe}, p.41, 4.2.~Theorem in particular), $X = \PN{n}$.

Let now $l$ be a line in $X = \PN{n}$. Decompose $E|_l = \oplus_{i=1}^r \O_l(a_i)$ and consider the natural inclusion $J_{k-1}(\O_{\PN{1}}(a_i)) \hookrightarrow J_1(J_1(\dots J_1(\O_{\PN{1}}(a_i))\dots))$. The bundle on the right can be computed from $J_1(\O_{\PN{1}}(l)) = \O_{\PN{1}}(l-1)^{\oplus 2}$ for $l \not= 0$ and $J_1(\O_{\PN{1}}) = \O_{\PN{1}} \oplus \O_{\PN{1}}(-2)$. This gives that all $a_i \ge k-1$, since $k \ge 2$ and we know from above, that $E$ is ample. \Cor{Chern} implies then $E|_l = \oplus_r \O_{\PN{1}}(k-1)$. This means the bundle $F = E \otimes \O_{\PN{n}}(-k+1)$ is trivial on every line in $X = \PN{n}$, hence trivial (this fact on $\PN{n}$ can be easily shown by induction on $n$.).
\end{proof2}

\

Note that the rest of the Theorem from the introduction now follows immediately using the numerical equivalence $\det E^* \equiv \frac{r(k-1)}{n+1}\cdot K_X$. Indeed, if we assume that $X$ is projective and that $\det E^*$ is not nef, then $K_X$ is not nef since $k > 1$, and $X$ contains a rational curve by the cone theorem (see for example \cite{MiPe}, p.37, 2.6.~Theorem). \Prop{RatCurve} implies $X \simeq \PN{n}$ and $E \simeq \oplus_r \O_{\PN{n}}(k-1)$ proving 1.) of the Theorem. Concerning 2.), if $\det E^* \equiv 0$, then $K_X \equiv 0$, $X$ admits a holomorphic affine connection and is covered by a complex torus. Finally, if $\det E^*$ is ample, then $K_X$ is ample. By the theorem of Aubin and Yau, $X$ admits a K\"ahler--Einstein metric so that $X$ is covered by the unit ball by \Theo{KE}.

\begin{corollary} \label{OX}
Let $X$ be a compact K\"ahler manifold. The first jet sequence of $\O_X$ is always split exact. One of the higher jet sequences
 \[0 \lra S^k\Omega^1_X \lra J_k(\O_X) \lra J_{k-1}(\O_X) \lra 0,\]
$k > 1$, splits if and only if $X$ is covered by a torus, i.e., admits a flat holomorphic affine connection. The splitting of one higher jet sequence of $\O_X$ implies the splitting of all jet sequences of $\O_X$.
\end{corollary}

\begin{proof}
By the main Theorem just proved the splitting of one of the higher jet sequences of $\O_X$ implies that $X$ is an \'etale quotient of a torus $T$. If $\mu: T \to X$ is \'etale, then $\mu^*J_i(\O_X) = J_i(\O_T)$ implying that the $i$--th jet sequence of $\O_X$ splits if and only if the $i$--th jet sequence of $\O_T$ is split exact. Since $\O_T$ is a theta characteristic on a torus, the result now follows from Biswas' result.
\end{proof}

\

For $X$ a curve, a K\"ahler surface or a projective threefold, \Cor{OX} has the following more general form: assume that there exists a theta characteristic $\theta \in \Pic(X)$. If one of the sequences
 \[0 \lra S^k \Omega_X^1 \otimes \theta^{1-k} \lra J_k(\theta^{1-k}) \lra J_{k-1}(\theta^{1-k}) \lra 0\]
splits for $k > 1$, then they all split and $X$ has a projective structure.


\end{document}